\documentclass[12pt]{amsart}
\usepackage{epsfig,amscd}
\usepackage[headings]{fullpage}

\newcommand{\psdraw}[2]
         {\begin{array}{c} \hspace{-1.3mm}
         \raisebox{-4pt}{\psfig{figure=#1.eps,width=#2}}
         \hspace{-1.9mm}\end{array}}

\theoremstyle{plain}

\newtheorem{theorem}{Theorem}
\newtheorem{proposition}{Proposition}[section]
\newtheorem{lemma}[proposition]{Lemma}
\newtheorem{corollary}[proposition]{Corollary}

\newtheorem{conjecture}{Conjecture}

\theoremstyle{definition}

\newtheorem{question}{Question}

\theoremstyle{remark}

\newtheorem{remark}[proposition]{Remark}

\def\printname#1{
    \if\draft y
        \smash{\makebox[0pt]{\hspace{-0.5in}
            \raisebox{8pt}{\tt\tiny #1}}}
    \fi
}

\newcommand{\cA}{{\mathcal A}}
\newcommand{\cP}{{\mathcal P}}
\newcommand{\cS}{{\mathcal S}}
\newcommand{\cT}{{\mathcal T}}

\newcommand{\cR}{{\mathcal R}}
\newcommand{\cX}{{\chi}}
\newcommand{\bR}{{\Bbb R}}
\newcommand{\bC}{{\mathbb C}}
\newcommand{\bZ}{{\Bbb Z}}

\newcommand{\tr}{\operatorname{tr}}

\def\gb{\mathfrak b}
\def\gt{\mathfrak t}
\def\gp{\mathfrak p}

\def\BZ{\mathbb Z}

\def\BR{\mathbb R}
\def\BC{\mathbb C}

\def\B{\mathcal B}

\def\cO{\mathcal O}

\def\R{\mathcal R}

\def\bb{\beta}

\def\S{\Sigma}

\def\la{\langle}
\def\ra{\rangle}

\def\b{\beta}

\def\sub{\subset}

\def\longto{\longrightarrow}

\def\pt{\partial}

\def\SL{\mathrm{SL}}

\newcommand{\lcr}{\raisebox{-5pt}{\mbox{}\hspace{1pt}
                  \epsfig{file=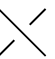}\hspace{1pt}\mbox{}}}
\newcommand{\ift}{\raisebox{-5pt}{\mbox{}\hspace{1pt}
                  \epsfig{file=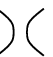}\hspace{1pt}\mbox{}}}
\newcommand{\zer}{\raisebox{-5pt}{\mbox{}\hspace{1pt}
                  \epsfig{file=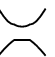}\hspace{1pt}\mbox{}}}

\usepackage{subfigure}
\begin{document}

\title[Jones Polynomial and A-polynomial]
      {The Colored Jones Polynomial and the A-Polynomial of Knots}

\author{Thang~T.~Q.~L\^e}
\address{
  School of Mathematics\\
  Georgia Institute of Technology\\
Atlanta, GA 30332-0160, USA } \email{letu@math.gatech.edu}
 \urladdr{http://www.math.gatech.edu/\~{}letu}

\thanks{
  This preprint is available electronically at
  {\tt http://www.math.gatech.edu/\~{}letu},
  and at {\tt http://xxx.lanl.gov/abs/math/0407521}.
}
\thanks{ The author  was
partially supported by an NSF grant}

\begin{abstract} We study relationships between the colored Jones
polynomial and the A-polynomial of a knot. The AJ conjecture (of
Garoufalidis) that relates the colored Jones polynomial and the
A-polynomial is established for a large class of two-bridge knots,
including all twist knots. We formulate a weaker conjecture and
prove that it holds for all two-bridge knots. Along the way we
also calculate the Kauffman bracket skein module of the
complements of two-bridge knots. Some properties of the colored
Jones polynomial are established.
\end{abstract}

\maketitle

\addtocounter{section}{-1}

\section{Introduction}

The Jones polynomial was discovered by Jones in 1984 \cite{Jones}
and has made a revolution in knot theory. Despite many efforts
little is known about the relationship between the Jones
polynomial and classical topology invariants like the fundamental
group. The $A$-polynomial of a knot, introduced in \cite{CCGLS},
describes more or less the representation space of the knot group
into $SL(2,\BC)$, and has been fundamental in geometric topology.

In the present paper we study relationships between the Jones
polynomial and the A-polynomial.  One main goal of the paper is to
establish for a large class of two-bridge knots the AJ conjecture
(made by Garoufalidis) that relates the colored Jones polynomial
and the A-polynomial. This class of knots contains for example all
twist knots, and much more. Another main result is the calculation
of the Kauffman bracket skein module of all two-bridge knots. This
generalizes the work \cite{Bullock,BullockL} where the
calculations were carried out for $(2p+1,2)$-torus knots and
twists knots, a special class of two-bridge knots. Our method is
more geometric  and we hope that it can be generalized to all
knots.

In a previous paper \cite{GL} Garoufalidis and the author proved
that for every knot, the colored Jones polynomial satisfies a
recurrence relation. The AJ conjecture states that when reducing
the quantum parameter to 1, the recurrence polynomial is
essentially equal to the $A$-polynomial (for details see below).
The present paper is independent of \cite{GL}, since we will prove
the existence of recurrence relations for two-bridge knots in
another way. We also formulate a weaker version of the AJ
conjecture (see Conjecture \ref{weakconjecture}) that we believe
reflects more accurately the algebra/topology relations between
the Jones polynomial and the A-polynomial. We prove that the
weaker conjecture holds true for all two-bridge knots.

Some properties  of the colored Jones polynomial are established.
We also show that for an arbitrary alternating knot, the degree of
the recurrence polynomial must be at least 2.

\subsection{The colored Jones polynomial and its recurrence ideal}

\subsubsection{The colored Jones polynomial}
For a knot $K$ in the 3-space $\BR^3\sub S^3$ the colored Jones
function (see for example \cite{MM,Tu})

$$J_K: \BZ \to  \cR := \BC[t^{\pm 1}]$$

is defined for integers $n\in \BZ$; its value $J_K(n)$ is known as
the colored Jones polynomial of the knot $K$ with color $n$. We
will recall the definition of $J_K(n)$ in section
\ref{jonescolor}.

In our joint work in S. Garoufalidis \cite{GL} we showed that the
function $J_K$ always satisfies a non-trivial recurrence relation
as described in the next subsection. Partial results were obtained
earlier by Frohman, Gelca, and Lofaro through their theory of
non-commutative A-ideal \cite{FGL,Gelca}, which also plays an
important role in the present paper.

\subsubsection{Recurrence  relations and $q$-holonomicity}
Consider  a function with domain the set of integers, $f: \BZ
\longto \cR$, and define the linear operators $L$ and $M$ acting
on such functions by:
$$
(M f)(n)= t^{2n}f(n), \hspace{1cm} (L f)(n) = f(n+1).
$$
It is easy to see that $LM= t^2ML$, and that $L^{\pm1},M^{\pm 1}$
generate the {\em quantum torus} $\cT$, a non-commutative ring
with presentation
$$
\cT=\cR\la M^{\pm 1},L^{\pm 1} \ra/(LM= t^2ML).
$$

We also use the notation $\cT_+$ for the subring of $\cT$ which
consists only of polynomials with non-negative powers of $M$ and
$L$. Traditionally $\cT_+$ is called the {\em quantum plane}.

The {\em recurrence ideal} of the discrete function $f$ is the
left ideal $\cA$ in $\cT$ that annihilates $f$:
$$
\cA=\{P \in \cT \, | \quad P\,f=0\}.
$$

We say that $f$ is $q$-{\em holonomic}, or $f$ satisfies a linear
recurrence relation, if $\cA\neq 0$. In \cite{GL} we proved that
for every knot $K$, the function $J_K$ is $q$-holonomic. Denote by
$\cA_K$ the recurrence ideal of $J_K$.

\subsubsection{An example} For the right-handed trefoil, one has

$$J_K(n) =
\frac{(-1)^{n-1}\, t^{2-2n}}{1-t^{-4}}\,\sum_{k=0}^{n-1}
t^{-4nk}\prod_{i=0}^k (1-t^{4i-4n}).$$

The function $J_K$ satisfies $p\, J_K =0$, where

$$ P = (t^4 M^{10} - M^6)L^2 - (t^2 M^{10}+ t^{-18}-
t^{-10} M^6 - t^{-14} M^4)L + (t^{-16} -t^{-4}M^4).$$

Together with the initial conditions $J_K(0)=0, J_K(1)=1$, this
recurrence relation determines $J_K(n)$ uniquely.

\subsubsection{Generator of the recurrence ideal}  The quantum torus $\cT$ is not a principal ideal
domain, and $\cA_K$ might not be generated by a single element.
Garoufalidis \cite{Ga2} noticed that by adding to $\cT$ all the
inverses of polynomials in $M$ one gets a principal ideal domain
$\tilde\cT$, and hence from the ideal $\cA_K$ one can define a
polynomial invariant. Formally one can proceed as follows. Let
$\cR(M)$ be the fractional field of the polynomial ring $\cR[M]$.
Let $\tilde \cT$ be the set of all Laurent polynomials in the
variable $L$ with coefficients in $\cR(M)$:

$$\tilde
\cT =\{\sum_{k\in \BZ}a_k(M) L^k \,\, | \quad a_k(M)\in \cR(M),
\,\,\, a_k=0  \quad \text{almost everywhere} \},
$$
and define the product in  $\tilde \cT$ by $a(M) L^k \cdot b(M)
L^l=a(M)\, b(t^{2k}M) L^{k+l}.$

Then it is known that every left ideal in $\tilde\cT$ is
principal, and $\cT$ embeds as a subring of $\tilde \cT$. The
extension $\tilde \cA_K  := \tilde \cT  \cA_K$  of $\cA_K$ in
$\tilde \cT$ is then generated by a single polynomial

$$\alpha_K(t;M,L) = \sum_{i=0}^n \alpha_{K,i}(t;M) L^i \,\in\, \cT_+,$$
where the degree in $L$ is assumed to be minimal and all the
coefficients $\alpha_{K,i}(t;M)\in \BZ[t^{\pm1},M]$ are assumed to
be co-prime. That $\alpha_K$ can be chosen to have integer
coefficients follows from the fact that $J_K(n) \in
\BZ[t^{\pm1}]$. It is clear that $\alpha_K(t;M,L)$ annihilates
$J_K$, and hence it is in the recurrence ideal $\cA_{K}$. Note
that $\alpha_K(t;M,L)$ is defined up to a factor $\pm t^a M^{b},
a,b\in \BZ$. We will call $\alpha_K$ the {\em recurrence
polynomial} of $K$. For example, the polynomial $P$ in the
previous subsection is the recurrence polynomial of the
right-handed trefoil.

\begin{remark} If $P$ is a polynomial in $t$ and $M$ (no
$L$), and $Pf=0$ then $P=0$. Hence adding all the inverses of
polynomials in $M$ does not affect the recurrence relations.
\end{remark}

\subsection{Main results} Let $\epsilon$ be the map reducing
$t=-1$. Formally, if $V$ is an $\cR$-module, then let $\epsilon
(V) = \BC \otimes_{\cR} V$, where $\BC$ is considered as an
$\cR$-module by setting $t=-1$. Also if $x\in V$ then
$\epsilon(x)$ is the image of $1\otimes x$ in $\epsilon(V)$. Thus
$\epsilon(\alpha_K)$ is the polynomial obtained from
$\alpha_K(t;M,L)$ by putting $t=-1$. For example, when $K$ is the
right-handed trefoil, $\epsilon(\alpha_K) = (M^4-1)(L-1)(LM^6+1)$.

For non-zero $f,g\in \BC[M,L]$, we say that $f$ is {\em
$M$-essentially equal to } $g$, and write $$f \,\overset{M}{=}\,
g,$$ if the quotient $f/g$ does not depend on $L$. We say that two
algebraic subsets of $\BC^2$ with parameters $(M,L)$  are {\em
$M$-essentially equal} if they are the same up to adding some
lines parallel to the $L$-axis. It's clear that if $f$ is
$M$-essentially equal to $g$, then $\{f=0\}$ and $\{g=0\}$ are
$M$-essentially equal. Here $\{f=0\}$ is the algebraic set of zero
points of $f$.

\subsubsection{The AJ Conjecture} Let $A_K\in \BZ[L,M]$ be the
$A$-polynomial of $K$ (see \cite{CCGLS,CL}); we will review its
definition in section \ref{Apolynomial}.    Garoufalidis
\cite{Ga2} made the following conjecture.

\begin{conjecture} {\rm ({\bf The AJ conjecture})}
 The polynomials $\epsilon(\alpha_K)$ and $(L-1)A_K$ are $M$-essentially equal.
\end{conjecture}

Actually, this is the strong version. The weak version of the
conjecture says that $\{\epsilon(\alpha_K)=0\}$ and
$\{(L-1)A_K=0\}$ are $M$-essentially equal. The algebraic set
$\{(L-1)A_K=0\}$ is known as the {\em deformation variety} of the
knot group, with the component $\{L-1=0\}$ corresponding to
abelian representations of the knot group into $SL_2(\BC)$, and
$\{A_K=0\}$ -- to non-abelian ones.

 Garoufalidis \cite{Ga2}
verified the conjecture for the trefoil and the figure 8 knot.
Takata \cite{Takata} gave some evidence to support the conjecture
for twist knots, but did not prove it. Both works are based
heavily on the computer programs of Wilf and Zeilberger. Hikami
\cite{Hikami} verifies the conjecture for torus knots. In all
these works direct calculations with explicit formulas are used.

In the present paper we prove the conjecture for a large class of
two-bridge knots, using a more conceptual approach. Two-bridge
knots $\gb(p,m)$ are parametrized by a pair of odd positive
integers $m<p$, with $\gb(p,m)= \gb(p,m')$ if $mm'=1 \pmod{p}$
(see \cite{Burde} and section \ref{twobridge} below).

\begin{theorem} Suppose $K=\gb(p,m)$ is a two-bridge knot.

a) The recurrence polynomial $\alpha_K$ has $L$-degree less than
or equal to $(p+1)/2$.

b) The algebraic set $\{\epsilon(\alpha_K)=0\}$ is $M$-essentially
equal to an algebraic subset of $\{ (L-1)A_K=0\}$.

c) The AJ conjecture holds true if

(*) the $A$-polynomial is $\BZ$-irreducible and has $L$-degree
$(p-1)/2$. \label{main1}
\end{theorem}

Here  $\BZ$-irreducibility means irreducibility in $\BZ[M,L]$. There
are many two-bridge knots that satisfy condition (*). For example in
a recent work \cite{Hoste} Hoste and Shanahan proved that  all the
twist knots satisfy the condition (*). Hence we have the following
corollary.

\begin{corollary}The AJ conjecture holds true for twist knots.
\end{corollary}

In a separate paper \cite{Le10} we will prove that if both $p$ and
$(p-1)/2$ are prime, then $\gb(p,m)$ satisfies the condition (*).
Also knot tables show that many two-bridge knots with small $p,m$
satisfy the condition (*).

\subsubsection{The Kauffman Bracket Skein Module of knot
complements}

Our proof of the main theorem is more or less based on the
ideology that the Kauffman Bracket Skein Module (KBSM) is a
quantization of the $SL_2(\BC)$-character variety (see
\cite{Bullock4,PS} and section \ref{peripheral} below), which has
been exploited in the work of Frohman, Gelca, and Lofaro
\cite{FGL} where they defined the non-commutative $A$-ideal. The
calculation of the KBSM of a knot complement is a difficult task.
Bullock \cite{Bullock} and recently Bullock and Lofaro \cite
{BullockL} calculated the KBSM for the complements of $(2,2p+1)$
torus knots and twist knots. Another main result of this paper is
a generalization of these works: We calculate explicitly the  KBSM
for complements of all two-bridge knots. We will use another, more
geometric approach that allows us to get the results for all
two-bridge knots.

\subsubsection{Other results} We also prove that the growth of degree (or breadth) of the
colored Jones polynomial of an arbitrary knot is at most quadratic
with respect to the color, and if the knot is alternating, then
the growth is exactly quadratic, given by explicit formula. This
is based on the exact estimate of the crossing number, used in the
proof (of Kauffman, Murasugi, and Thistlethwaite) of the Tait
conjecture on the crossing number of alternating knots. As a
corollary, we show that the $L$-degree of the recurrence
polynomial of an alternating knot must be at least 2.

\subsection{Plan of the paper} In Section 1 we
review the theory of skein modules the colored Jones polynomial. In
Section 2 we study the growth of the degree of the colored Jones
polynomial and the $L$-degree of the recurrence polynomial. In
Section 3 we review the A-polynomial and  introduce a closely
related polynomial, $B_K$. Section 4 is devoted to the ``quantum"
version of $B_K$, the peripheral polynomial. We will formulate
another weaker version of the AJ conjecture and prove it holds true
for two-bridge knots (in Section \ref{proof1}). In Section 5 we
calculate the skein module of the complement of two-bridge knots.
The last section contains a proof of the Theorem \ref{main1}.

\subsection{Acknowledgments} I would like to thank S. Garoufalidis,
R. ~Gelca, T. ~Ohtsuki, A. ~Referee, P. ~Shanahan, A. ~Sikora, T.
Takata, and X. ~Zhang for valuable helps and discussions.

\section{The colored Jones polynomial and skein modules}
\label{jonescolor}
 We recall the definition and known facts
about the colored Jones polynomial through the theory of Kauffman
Bracket Skein Modules which was introduced by Przytycki and
Turaev, see the survey \cite{Przytycki}.

\label{Preliminaries}

\subsection{Skein modules}

 Recall that $\cR=\bC[t^{\pm1}]$.
 A {\em framed link} in an oriented  3-manifold $Y$ is  a disjoint union of embedded circles, equipped with a
 non-zero normal vector field. Framed links are considered up to isotopy. In all figures we will draw  framed links,
or part of them, by lines as usual, with the convention that the
framing is blackboard. Let $\mathcal{L}$ be the set of isotopy
classes of framed links in the manifold $Y$, including the empty
link. Consider the free $\cR$-module with basis $\mathcal{L}$, and
factor it by the smallest submodule containing all expressions of
the form $\displaystyle{\lcr-t\zer-t^{-1}\ift}$ and
$\bigcirc+(t^2+t^{-2}) \emptyset$, where the links in each
expression are identical except in a ball in which they look like
depicted. This quotient is denoted by $\cS(Y)$ and is called the
Kauffman bracket skein module, or just skein module, of  $Y$.

When $Y= \Sigma \times [0,1]$, the cylinder over the surface
$\Sigma$, we also use the notation $\cS(\S)$ for $\cS(Y)$. In this
case $\cS(\S)$ has an algebra structure induced by the operation
of gluing one cylinder on top of the other. The operation of
gluing the cylinder over $\partial Y$ to $Y$ induces a
$\cS(\partial Y)$-left module structure on $\cS(M)$.

\subsection{Example: $\cS(\bR^3)$ and the Jones polynomial}

When $Y$ is  the 3-space $\bR^3$ or the 3-sphere $S^3$, the skein
module $\cS(Y)$ is free over $\cR$ of rank one, and is spanned by
the empty link. Thus if $\ell$ is a framed link in $\bR^3$, then
its value in the skein module $\cS(\bR^3)$ is $\langle \ell
\rangle$ times the empty link, where $\langle \ell \rangle \in
\cR$, known as the Kauffman bracket of $\ell$ (see
\cite{Kauffman,Lickorish}), and is just the Jones polynomial of
{\em framed links} in a suitable normalization.
\subsection{Example: The solid torus and the colored Jones
polynomial}

The solid torus $ST$ is the cylinder over an annulus, and hence
its skein module $\cS(ST)$ has an algebra structure. The algebra
$\cS(ST)$ is the polynomial algebra $\cR[z]$ in the variable $z$,
which is a knot representing the core of the solid torus.

Instead of the $\cR$-basis $\{1, z, z^2,\dots\}$, two other bases
are often useful. The first basis consists of the Chebyshev
polynomials $T_n(z)$, $n\geq 0$,
 defined by $T_0(x)=2$, $T_1(x)=x$, and $T_{n+1}(x)=xT_n-T_{n-1}$.
  The second basis consists of polynomials
$S_n(z)$, $n\geq 0$, satisfying the same recurrence relation, but
with $S_0(x)=1$ and $S_1(x)=x$.  Extend both polynomials by the
recurrence relation to all indices $n\in{\mathbb Z}$. Note that
$T_{-n}=T_n$, while $S_{-n}=-S_{n-2}$.

For a framed knot $K$ in a 3-manifold $Y$ we define the $n$-th
power $K ^n$ as the link consisting of $n$ parallel copies of $K$.
Using these powers of a knot,  $S_n(K)$ is defined as an element
of $\cS(Y)$. In particular, if $Y=\bR^3$ one can calculate the
bracket $\langle S_n(K)\rangle\in \cR$, and it is essentially the
colored Jones polynomial. More precisely, we will define the
colored Jones polynomial $J_K(n)$ by the equation

$$ J_K(n+1):=(-1)^{n} \times\langle S_n(K) \rangle.  $$

The $(-1)^n$ sign is added so that when $K$ is the trivial knot,
$$J_K(n) = [n] := \frac{t^{2n} -t^{-2n}}{t^2-t^{-2}}.$$ Then
$J_K(1)=1$, $J_K(2)= - \langle K \rangle$. We extend the
definition for all integers $n$ by $J_K(-n)= -J_K(n)$ and
$J_K(0)=0$. In the framework of quantum invariants, $J_K(n)$ is
the $sl_2$-quantum invariant of $K$ colored by the $n$-dimensional
simple representation of $sl_2$.

We will always assume $K$ has 0 framing. In this case $J_K(n)$
contains only even powers of $t$, i.e. $J_K(n) \in \BZ[t^{\pm2}]$.
Hence the recurrence polynomial $\alpha_K$ can be assumed to have
only even powers of $t$.

\subsection{Example: cylinder over the torus and the
non-commutative torus}

A pair of {\em oriented} meridian and longitude on the torus
${\mathbb T}^2$ will define an {\em algebra} isomorphism $\Phi$
between $\cS({\mathbb T}^2)$ and the symmetric part of the quantum
torus $\cT$ as follows.

For a pair of integers $a,b$ let $(a,b)_T = T_d((a',b')_T)$, where
$d$ is the greatest common divisor of $a$ and $b$, with $a=da',
b=db'$, and $(a',b')_T$ is the closed curve without
self-intersection on the torus that is homotopic to $a'$ times the
meridian plus $b'$ times the longitude. Here $T_d$ is the Chebyshev
polynomial defined above; and the framing of a curve on ${\mathbb
T}^2$ is supposed to be parallel to the surface ${\mathbb T}^2$.
Note that in the definition of skein modules we use {\em
non-oriented links}, hence $(a,b)_T= (-a,-b)_T$. As an $\cR$-module,
$\cS({\mathbb T}^2)$ is the quotient of the free $\cR$-module
spanned by $\{(a,b)_T, (a,b)\in \bZ^2\}$ modulo the relations
$(a,b)_T=(-a,-b)_T$.

Recall that the quantum torus $\cT$ is defined as $\cT =
\cR\langle L^{\pm 1},M^{\pm 1}\rangle /(LM=t^2ML)$. Let
$\cT^\sigma$ be the subalgebra of $\cT$ invariant under the
involution $\sigma$, where $\sigma(M^aL^b)=M^{-a}L^{-b}$. Frohman
and Gelca in \cite{FG} showed that the map

$$ \Phi: \cS({\mathbb T}^2) \to \cT^\sigma,
\qquad \Phi((a,b)_T) = (-1)^{a+b}t^{ab}\, (M^aL^b + M^{-a}
L^{-b})$$ is an isomorphism of {\em algebras}.

\label{qTorus}

\subsection{Two-punctured disk} \label{d3}
\begin{figure}[htpb]
$$ \psdraw{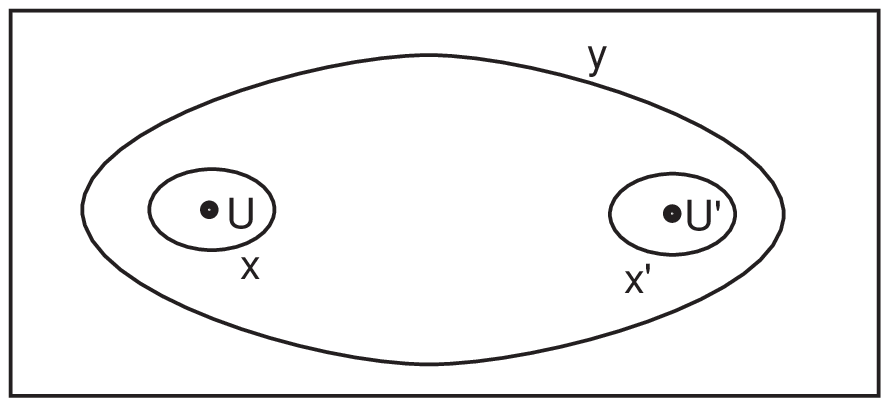}{2.5in} $$
\caption{The 2-punctured disk}\label{disk}
\end{figure}

Let $F$  be
 the rectangle (in $\bR^2$)
$$\{ (x_1,x_2) \in \BR^2 \mid 0\le x_1 \le 6, 0\le x_2 \le 2\},$$
without 2 interior points $U = (1,1)$ and $U'=(5,1)$. Then
$\cS(F)$ is the polynomial algebra $\cR[x,x',y]$, where $x$ is a
small loop around $U$, $x'$ a small loop around $U'$, and $y$ a
loop circling both $U,U'$ (see Figure \ref{disk}). Note that in
this case $\cS(F)$ is a commutative algebra, which is not true in
general when $F$ is replaced by an arbitrary surface, see
\cite{Przytycki}.

\begin{figure}[htpb]
$$ \psdraw{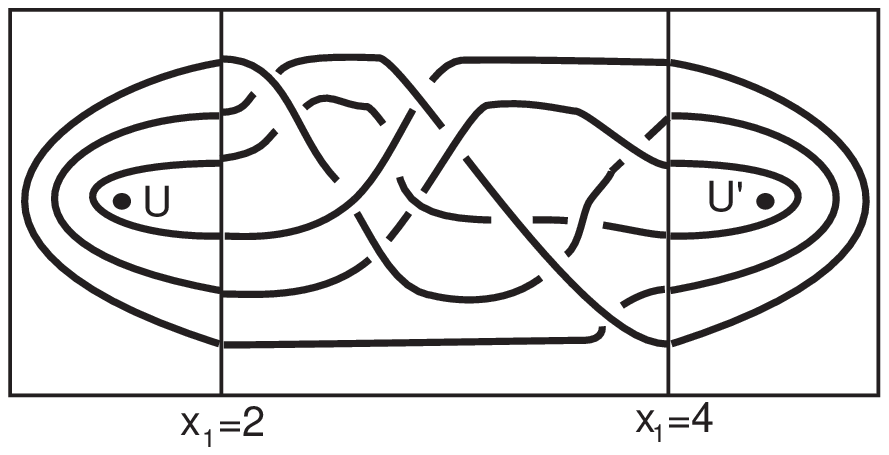}{3.3in} $$
\caption{The braid $\beta$ (in the middle) and its
closure}\label{closure}
\end{figure}

We will identify $F$ with the section $F\times {1/2}$ in $F\times
[0,1]$.
 The lines $x_1=2$ and $x_1=4$ divide $F$ into 3
parts: the left part $F_l$ (containing $U$), the right $F_r$
(containing $U'$), and the middle part $F_m$. Let $\beta$ be a braid
on $2k$ strands in $F_m\times [0,1]$ with boundary points on the
lines $x_1=2$ and $x_1=4$ (on $F$), i.e. $\beta$ consists of $2k$
connected paths, each begins at a point on $x_1=2$ (on $F$) and goes
monotonously to the right until ending at a point on $x_1=4$, see
Figure \ref{closure}.  Let the closure $\hat \beta$ be the link in
$F\times [0,1]$ obtained by closing $\beta$ using $k$ parallel
rainbow arcs in $F_l$ and $k$ parallel rainbow arcs in $F_r$, as
shown in Figure \ref{closure}.

\begin{lemma} As an
element of the skein module $\cS(F)= \cR[x,x',y]$ the closure link
$\hat \beta$ of a braid on $2k$ strands is a polynomial having
$y$-degree $k$, and the coefficient of $y^k$ is invertible and of
the form $\pm t^l, l\in \BZ$. \label{disk21}
\end{lemma}

\begin{proof} Consider the diagram $D$ of $\hat \beta$ on $F$. At every
crossing point there are 2 ways  to smooth the diagram, one
positive that gives coefficient $t$ in the skein relation, and the
other is negative. A {\em state} is the result of smoothing all
the crossing; what one has is a bunch of non-intersecting circles,
each is one of  $x,x',y$, or the trivial loop. The link $\ell$,
considered as an element of $\cS(F)\equiv \cR[x,x',y]$, is the sum
over all states, each with coefficient a power of $t$. It is clear
that for the closure $\hat \beta$ of a braid, there is only one
state (the one with horizontal resolutions everywhere) that gives
the maximal power $y^k$, and its coefficient is a power of $t$.
\end{proof}
\label{disk2}

\section{Some properties of the colored Jones polynomial}
\subsection{} One of the best known applications of the Jones
polynomial is a proof (Kauffman, Murasugi, and Thistlethwaite) of
the Tait conjecture on the crossing number of alternating links,
based on an exact estimate of the crossing number  using the {\em
breadth} of the Jones polynomial. We will need a generalization of
this estimate for the colored Jones polynomial.

For a Laurent polynomial $P(t)\in \BZ[t^{\pm1}]$ let $d_+(P)$ and
$d_-(P)$ be respectively the maximal and minimal degree of $t$ in
$P$. The difference $br:= d_+-d_-$ is called the {\em }breadth of
$P$. For a link diagram $D$, let $s_+(D)$ and $ s_-(D)$ be the
number of circles obtained by positively (respectively,
negatively) smoothing all the double points.

\begin{proposition}
a) Suppose $K$ is a knot with a knot diagram $D$ having $k$
crossings. Assuming the framing is black-board. Then

$$\begin{aligned} d_+(J_K(n)) &\le k(n-1)^2 + 2(n-1) s_+(D),\\
d_-(J_K(n)) &\ge -k(n-1)^2 - 2(n-1) s_-(D).
\end{aligned}$$
Hence the breadth of $J_K(n)$ grows at most as a quadratic
function in $n$.

 b) If $K$ is a non-trivial alternating knot with
$k$ crossings. Then the breadth of $J_K(n)\in \BZ[t^{\pm1}]$ is a
quadratic polynomial in $n$. Moreover
$$ br(J_K(n))=  2 k(n-1)^2 + 2(n-1)(k+2).$$
\label{degree10}
\end{proposition}
\begin{proof} a) The $n$-parallel $D^n$ of $D$
will have $kn^2$ double points. In addition, it is easy to see
that $s_\pm (D^n) = n s_\pm (D)$. Hence Lemma 5.4 of
\cite{Lickorish} says that

$$ d_+\langle D^n \rangle  \le  f(n):= k(n-1)^2 + 2 (n-1) s_+.$$
Note that here we have taken into account the fact that we shift
the index $n\to n+1$, and use the normalization in which the
unknot takes value $-t^2-t^{-2}$. Note that $f(n)$ is a strictly
increasing function, $f(n+1) > f(n)$. Recall that $S_n(K) = D^n +
\text{terms of lower degrees in $D$}$. Hence one has

$$ d_+\langle S_n(D) \rangle   \le k(n-1)^2 + 2 (n-1) s_+.$$
The proof for $d_-$ is similar.

b)
  For an alternating knot diagram one
has $s_++s_-=k+2$. Note that the diagram $D^n$ is adequate in the
sense of \cite[Chapter 5]{Lickorish}. For adequate diagrams the
estimates of $d_\pm$ in part a) are exact (see the proof of Lemma
5.4 of \cite{Lickorish}) and the result follows immediately.
\end{proof}

\subsection{The $L$-degree of the recurrence polynomial}

\begin{proposition}Suppose $K$ is a non-trivial alternating knot.
Then the recurrence polynomial $\alpha_K$ has $L$-degree greater
than 1. \label{deg2}
\end{proposition}
\begin{proof} Assume the contrary that
$ \alpha_K = P(t;M)\, L + P_0(t;M)$, where $P,P_0\in
\BZ[t^{\pm1},M^{\pm1}]$. Garoufalidis \cite{Ga2} showed that the
polynomial $\sigma(\alpha_K)= P(t;M^{-1})\, L^{-1} +
P_0(t;M^{-1})$ is also in the recurrence ideal. Since $\alpha_K$
is {\em the} generator, it follows that for some $\gamma(t;M) \in
\cR(M)$
$$L\,\sigma(\alpha_K) =
\gamma(t;M) \, \alpha_K.$$ One can then easily show that, after
normalizing both $P,P_0$  by a same power of $M$, one has
$$P_0(t;M)= P(t; t^{-2}
M^{-1}).$$

The equation $\alpha_K J_K=0$ can now be rewritten as

$$ J_K(n+1) = -\frac{P(t;t^{-2-2n})}{P(t;t^{2n})} J_K(n).$$

It is easy to see that for $n$ big enough, the difference of the
breadths $br(P(t;t^{-2-2n}))- br(P(t;t^{2n}))$ is a constant
depending only on the polynomial $P(t;M)$, but not on $n$. From
the above equation it follows that the breadth of $J_K(n)$, for
$n$ big enough, is a linear function on $n$. This contradicts
Proposition \ref{degree10}.
\end{proof}

\subsection{More on the recurrence polynomial of knots}
The following proposition was known to Garoufalidis.

\begin{proposition} When reduced by $t=-1$, the recurrence
polynomial $\alpha_K$ is divisible by $L-1$. In other words,
$\displaystyle{\frac{\epsilon(\alpha_K)}{L-1}}\in \BZ[M,L]$.
\label{divi}
\end{proposition}

\begin{proof} For a function $ f(t^2,n)$ of two variables $t^2 \in
\BC$ and $n\in \BZ$ let $\overline f(z)$ be the limit of
$f(t^2,n)$ when

($\dag$) $t^2 \to 1$ and $t^{2n}$ is kept equal to $z$ all the
time.

The Melvin-Morton conjecture \cite{MM}, proved by Bar-Natan and
Garoufalidis \cite{BG}, showed that $h(z):=\overline {J_K(n)}$ is
the inverse of the Alexander polynomial. In particular, $h(z)\not
= 0$.

Lemma \ref{new} below shows that $\overline{ J_K(n+k)} =
\overline{ J_K(n)}$ for any fixed $k$. Hence the operator $L$
becomes the identity after taking the limit ($\dag$). Thus
applying the limit ($\dag$) to the equation $ \alpha_K J_K=0$ we
see that

$$ \alpha_K|_{t^2=1, M=z, L=1} \, h(z)=0.$$

Since $h(z) \not = 0$, one has $\alpha_K|_{t^2=1, M=z, L=1}=0$,
which is equivalent to  the lemma.
\end{proof}

\begin{lemma} For every fixed integer $k$, one has $\overline{ J_K(n+k)} =
\overline{ J_K(n)}$. \label{new}
\end{lemma}

\begin{proof} For a knot $K$ of framing 0, $J_K(n)/[n]$ is a
Laurent polynomial in $t^4$, and  (see \cite{MM})

$$ \frac{J_K(n)}{[n]}|_{t^4= \exp\hbar}= \sum_{l=0}^\infty P_{l}(n) \hbar^l,$$
where $P_l(n)$ is a polynomial in $n$ of degree at most $l$:

$$P_l(n) = P_{l,l} n^l + P_{l,l-1} n^{l-1} +\dots P_{l,1} n +
P_{l,0}.$$

The limit ($\dag$) is the same as the limit $n\to \infty$, with
$\hbar = 2\ln z/n$. Under this limit,

 $$ (n+k)^i \, \hbar^k \to \begin{cases}
 0 & \text{if $i<l$}\\
 (2\ln z)^l & \text{if $i=l$}
 \end{cases},$$
 which does not depend on $k$. \end{proof}

\section{The A-polynomial and  its sibling}
\label{Apolynomial}

We briefly recall here the definition of the A-polynomial and
introduce a sibling of it. We will say that $f$ is {\em
$M$-essentially divisible by} $g$ if $f$ is $M$-essentially equal
to a polynomial divisible by $g$.

\subsection{The character variety of a group}

The set of representations of a finitely presented group $\pi$
into $\SL_2(\BC)$ is an algebraic set defined over $\BC$, on which
$\SL_2(\BC)$ acts by conjugation. The naive quotient space, i.e.
the set of orbits, does not have a good topology/geometry. Two
representations in the same orbit (i.e. conjugate) have the same
character, but the converse is not true in general. A better
quotient, the algebro-geometric quotient denoted by $\cX(\pi)$
(see \cite{Culler,Lubotzky}), has the structure of an algebraic
set. There is a bijection between $\cX(\pi)$ and the set of all
characters of representations of $\pi$ into $SL_2(\BC)$, hence
$\cX(\pi)$ is usually called the {\em character variety} of $\pi$.
For a manifold $Y$ we use $\cX(Y)$ also to denote $\cX(\pi_1(Y))$.

Suppose  $\pi=\BZ^2$, the free abelian group with 2 generators.
Every pair  of generators $\lambda,\mu$ will define an isomorphism
between $\cX(\pi)$ and $(\BC^*)^2/\tau$, where $(\BC^*)^2$ is the
set of non-zero complex pairs $(L,M)$ and $\tau$ is the involution
$\tau(M,L)=(M^{-1},L^{-1})$, as follows: Every representation is
conjugate to an upper diagonal one, with $L$ and $M$ being the
upper left entry of $\lambda$ and $\mu$ respectively. The
isomorphism does not change if one replaces $(\lambda,\mu)$ with
$(\lambda^{-1},\mu^{-1})$.

\subsection{The A-polynomial}  Let $X$ be the closure of $S^3$
minus a tubular neighborhood $N(K)$ of a knot $K$. The boundary of
$X$ is a torus whose fundamental group  is free abelian of rank
two. An orientation of $K$ will define a unique pair of an
oriented meridian  and an oriented longitude such that the linking
number between the longitude and the knot is 0. The pair provides
an identification of $\cX(\pi_1(\partial X))$ and $(\BC^*)^2/\tau$
which actually does not depend on the orientation of $K$.

The inclusion $\pt X \hookrightarrow X$ induces the restriction
map

$$\rho : \cX(X) \longto \cX(\pt X)\equiv (\BC^*)^2/\tau$$
 Let $Z$ be the image of
$\rho$ and  $\hat Z \sub (\BC^*)^2$ the lift of $Z$ under the
projection $(\BC^*)^2 \to (\BC^*)^2/\tau$. The Zariski closure of
$\hat Z\sub (\BC^*)^2 \sub \BC^2$ in $\BC^2$ is an algebraic set
consisting of components of dimension 0 or 1. The union of all the
1-dimension components is defined by a single polynomial $A'_K\in
\BZ[M,L]$, whose coefficients are co-prime. Note that $A'_K$ is
defined up to $\pm 1$. It is known that $A'_K$ is divisible by
$L-1$, hence $A'_K= (L-1)A_K$, where $A_K\in \BC[M,L]$ is called
the $A$-polynomial of $K$. It is known that $A_K\in \BZ[M^2,L]$.
By definition, $A_K$ does not have repeated factor, and is not
divisible by $L-1$.

\begin{question} Can $A_K(L,M)$ have a factor a non-constant polynomial
depending on $M$ only? \label{question}
\end{question}

\subsection{The dual picture} It's also instructive and
convenient to see the dual picture in the construction of the
A-polynomial. For an algebraic set $Y$ let $R[Y]$ denote the ring
of regular functions on $Y$.  For example, $R[(\BC^*)^2/\tau]=
\gt^\sigma$, the $\sigma$-invariant subspace of  $\gt:=\BC[L^{\pm
1},M^{\pm 1}]$, where $\sigma(M^aL^b):= M^{-a}L^{-b}$.

The map $\rho$ in the previous subsection induces an algebra
homomorphism

$$\theta: R[\cX(\pt X)]  \equiv \gt^\sigma \longrightarrow
R[\cX(X)].$$

We will call the kernel  $\gp$ of $\theta$ the {\em classical
peripheral ideal}; it is an ideal of $\gt^\sigma$. Let $\hat \gp:=
\gt \, \gp$ be the ideal extension of $\gp$ in $\gt$. The set of
zero points of $\hat \gp$ is the closure of $\hat Z$ in $\BC^2$.

\subsection{A sibling of the A-polynomial}

The ring $\gt= \BC[M^{\pm1},L^{\pm 1}]$ embeds naturally into the
principal ideal domain $\tilde \gt:=\BC(M)[L^{\pm1}]$, where
$\BC(M)$ is the fractional field of $\BC[M]$. The ideal extension
of $\hat \gp$ in $\tilde \gt$, which is $\tilde \gt \, \hat \gp =
\tilde \gt \, \gp$,  is thus generated by a single polynomial
$B_K\in \BZ[M,L]$ which has co-prime coefficients and is defined
up to a factor $\pm M^a$ with $a\in \BZ$. Again $B_K$ can be
chosen to have integer coefficients because everything can be
defined over $\BZ$.

From the definitions one has immediately

\begin{proposition} The polynomial $B_K$ is $M$-essentially
divisible by $A'_K=(L-1)A_K$. The two algebraic sets $\{B_K=0\}$
and $\{A'_K=0\}$ are $M$-essentially equal. \label{500}
\end{proposition}

Note that $B_K$ might not be $M$-essentially equal to $A'_K$
because $B_K$ might contain repeated  factors. If the answer to
Question \ref{question} is negative, then the two algebraic sets
$\{B_K=0\}$ and $\{A'_K=0\}$ are equal, and we have a natural way
to define the multiplcity of factors of the $A$-polynomial, using
the $B$-polynomial.

\section{The quantum peripheral ideal and the peripheral polynomial}

\subsection{Skein modules as quantum deformations of character
varieties} Recall that $\epsilon$ is the map reducing $t=-1$. One
important result (Bullock, Przytycki, and Sikora
\cite{Bullock,PS}) in the theory of skein modules is that
$\epsilon (\cS(Y))$ has a natural algebra structure and, when
factored by its nilradical, is canonically isomorphic to
$R[\cX(Y)]$, the ring of regular functions on the character
variety of $\pi_1(Y)$. The product of 2 links in $\epsilon
(\cS(Y))$ is their union. Using the skein relation with $t=-1$, it
is easy to see that the product is well-defined, and that the
value of a knot in the skein module depends only on the homotopy
class of the knot in $Y$. The isomorphism between
$\epsilon(\cX(Y))$ and $R[\cX(Y)]$ is given by $K(r)= -\tr r(K)$,
where $K$ is a homotopy class of a knot in $Y$, represented by an
element, also denoted by $K$, of $\pi_1(Y)$, and $r:\pi_1(Y) \to
SL_2(\BC)$ is a representation of $\pi_1(Y)$.

In many cases the nilradical of $\epsilon (\cS(Y))$ is trivial,
and hence $\epsilon (\cS(Y))$ is exactly equal to the ring of
regular functions on the character variety of $\pi_1(Y)$. For
example, this is the case when $Y$ is a torus, or when $Y$ is the
complement of a two-bridge knots (see section \ref{twobridge}).

In light of this fact, one can consider $\cS(Y)$ as a quantization
of the character variety. \label{peripheral}

\subsection{The quantum peripheral ideal} Recall that  $X$ is the
closure of the complement of a tubular neighborhood $N(K)$  in
$S^3$. The boundary $\pt X$ is a torus, and using the preferred
meridian and longitude we will identify $\cS(\pt X)$ with
$\cT^\sigma$, see subsection \ref{qTorus}.

The embedding of $\partial X$ into $X$ gives us a map $\Theta:
\cS(\partial X)\equiv \cT^\sigma \longrightarrow \cS(X)$, which
can be considered as a quantum analog of $\theta$. One has the
following commutative diagram

$$ \begin{CD} \cT^\sigma  @>\Theta >> \cS(X) \\
@V \epsilon VV  @V \epsilon VV \\
\gt^\sigma  @>\theta >> R[\cX(X)]
\end{CD} $$

 The kernel of $\Theta$, denoted by $\cP$,  is called the {\em quantum peripheral ideal}; it is a left ideal
of $\cT^\sigma$ and can be considered as a quantum analog of the
classical peripheral ideal $\gp = \ker \theta$. The ideal $\cP$
was introduced by Frohman, Gelca, and Lofaro in \cite{FGL} and
there it is called simply the peripheral ideal. From the
commutative diagram it is clear that $\epsilon(\cP) \subset \gp$.
The following question is important.
\begin{question}
 Is it true that $\epsilon (\cP)$  can never be 0?
 \end{question}

\subsection{The peripheral polynomial}

Let us adapt the construction of the $B_K$ polynomial to the
quantum setting. Recall that $\tilde \cT$ (see Introduction) is a
principal left-ideal domain that contains $\cT^\sigma$ as a
subring. The left-ideal extension $\tilde \cP :=\tilde \cT \cP$ in
$\tilde \cT$ is generated  by a polynomial

$$\beta_K(t;M,L) = \sum_{i=0}^s \beta_{K,i}(t,M) \, L^i \in \cT_+,$$
where $s$ is assumed to be minimum and all the coefficients
$\beta_{K,i}(t,M) \in \BZ[t^{\pm 1}, M^{\pm 1}]$ are co-prime. We
call $\beta_K$ the {\em peripheral polynomial} of $K$, which is
defined up to $\pm t^a M^b$ with $a, b \in \BZ$.

\begin{proposition} $\epsilon(\beta_K)$ is $M$-essentially
divisible by $B_K$, and hence is $M$-essentially divisible by
$A'_K= (L-1)A_K$. \label{323}
\end{proposition}

\begin{proof} The
proposition follows the fact that $\epsilon \cP \sub \gp$ and
Proposition \ref{500}.
\end{proof}

\begin{proposition} Suppose $\epsilon(\cP)= \gp$. Then
$\epsilon(\beta_K) \overset{M}{=} B_K$. \label{t1}
\end{proposition}

\begin{proof}For the extensions $\hat \cP:= \cT \cP$ and $\hat \gp:= \gt
\gp$ we also have $\epsilon (\hat \cP) = \hat \gp$, since
$\epsilon \hat \cP = \epsilon (\cT \cP)= \epsilon(\cT)
\epsilon(\cP)=\gt\gp$.

From the definition we have that $h(M)\B_K \in \hat \gp$ for some
non-zero polynomial $h(M)\in \BZ[M]$. Hence $\epsilon^{-1}
(h(M)B_K) \sub \hat \cP$ is not empty. Take an element $u\in
\epsilon^{-1} (h(M)B_K)$; it is $M$-essentially divisible by
$\beta_K$, the generator. Applying the map $\epsilon$, one gets
$B_K$ is $M$-essentially divisible by $\epsilon(\beta_K)$.
Combining with Proposition \ref{323} one has $\epsilon(\beta_K)
\overset{M}{=} B_K$.
\end{proof}

\begin{conjecture} For every knot  we have $\epsilon(\cP)= \gp$,
and hence $\epsilon(\beta_K) \overset{M}{=} B_K$.
\label{weakconjecture}
\end{conjecture}

Later  we will show that for all two-bridge knots holds true the
conjecture, which is closely related to the AJ conjecture.  A
sufficient condition for the conjecture to hold true is given in
Section \ref{proof1}.

\subsection{The orthogonal ideal and recurrence relations} There
is a bilinear pairing

\begin{equation} \cS(N(K)) \otimes \cS(X) \to \cR, \quad \text{with}\quad  \ell \otimes \ell' \to \langle \ell ,\ell' \rangle := \langle \ell \cup \ell' \rangle
\in \cS(S^3)= \cR,\label{pair} \end{equation}

where $\ell$ and $\ell'$ are framed links in $N(K)$ and  $X$
respectively. The {\em orthogonal ideal} $\cO$ is defined by
$$ \cO := \{ \ell'\in \cS(\partial X)\quad  \mid \quad  \langle \ell ,\Theta(\ell')
\rangle=0 \quad \text{for every }\ell \in \cS(N(K))\}.$$

It's clear that $\cO$  is a left ideal of $\cS(\partial X)\equiv
\cT^\sigma$ and  $\cP \subset \cO$. In \cite{FGL}, where $\cO$ was
first introduced, $\cO$ was called the formal ideal. What is
important for us is the following
\begin{proposition} The orthogonal ideal is in the recurrence ideal of a
knot, $\cO \sub \cA_K$. As a consequence, $\cP \sub \cA_K$.
\label{inclu}
\end{proposition}

This was proved by Garoufalidis \cite{Ga2}. Frohman, Gelca, and
Lofaro \cite{FGL}  proved that every element $\ell'$ in the
orthogonal ideal $\cO$ gives rise to a linear recurrence relation
for the colored Jones polynomial. The idea is simple and
beautiful: $\ell'$ annihilates everything in $\cS(N(K))$, in
particular, $\langle T_n(z),\Theta(\ell')\rangle =0$; but this
equation, after some calculation, can be rewritten as a linear
recurrence relation for the colored Jones polynomial.
Garoufalidis, using the Weyl symmetry, further simplified the
recurrence relation, and obtained  that $\cO= \cA_{K}\cap
\cT^\sigma$, which is stronger than the proposition.

\begin{conjecture} The right kernel of the bilinear form
(\ref{pair}) is trivial.
\end{conjecture}
This conjecture implies that $\cO=\cP$, from which, due to an
argument of A. Sikora and the author, one can show that the colored
Jones polynomial distinguishes the unknot from other knots.

\subsection{Relation between the peripheral and recurrence
polynomials}

\begin{lemma}
The peripheral polynomial $\beta_K$ is divisible by the recurrence
polynomial $\alpha_K$ in the sense that there are polynomials
$g(t,M) \in \BZ[t,M]$ and $\gamma(t,M,L) \in \cT_+$ such that

\begin{equation} \beta_K(t,M,L) = \frac{1}{g(t,M)}\, \gamma(t,M,L) \,
\alpha_K(t,M,L).\label{sao}\end{equation}

Moreover $g(t,M)$ and $\gamma(t,M,L)$ can be chosen so that
$\epsilon g\not = 0$.
\end{lemma}

\begin{proof} From Proposition \ref{inclu} we have that $\cP \sub
\cA$. Hence the left-ideal extension $\hat \cP:=\cT \cP$ is also a
subset of $\cA$, since both are left ideals of $\cT$. It follows
that $\beta_K$, as the generator of the extension of $\hat\cP$ in
$\tilde\cT$, is divisible by the generator of the extension of
$\cA$, and (\ref{sao}) follows.

We can assume that $t+1$ does not divide both $g(t,M)$ and
$\gamma(t,M,L)$ simultaneously. If $\epsilon g = 0$ then $g$ is
divisible by $t+1$, and hence $\gamma$ is not. But then from the
equality

$$ g \beta_K = \gamma \, \alpha_K,$$
it follows that $\alpha_K$ is divisible by $t+1$, which is
impossible, since all the coefficients of powers of $L$ in
$\alpha_K$ are supposed to be co-prime.
\end{proof}

\section{Two-bridge knots and their skein modules}

\label{twobridge}

\newcommand{\cc}{\mathfrak{u}}
\newcommand{\cb}{\mathfrak{v}}
\newcommand{\ft}{\omega}
\newcommand{\fv}{\mathfrak{n}}
\newcommand{\cK}{\mathcal{K}}
\newcommand{\fk}{\mathfrak{a}}

\subsection{Two-Bridge knots}

 A two-bridge knot is a knot $K\subset S^3$ such that
there is a 2-sphere $S^2\subset S^3$ that separates $S^3$ into 2
balls $B_1$ and $B_2$, and the intersection of $K$ and each ball
is isotopic to 2 trivial arcs in the ball. The branched double
covering of $S^3$ along a two-bridge knot is a lens space
$L(p,m)$, which is obtained by doing a $p/m$ surgery on the
unknot. Such a two-bridge knot is denoted by $\gb(p,m)$. It is
known that both $p,m$ are odd. One can always assume that $ p> m
\ge 1$. It is known that $\gb(p,m) = \gb(p,m')$ if $m m' \equiv 1
\pmod{p}$.

\begin{figure}[htpb] $$ \psdraw{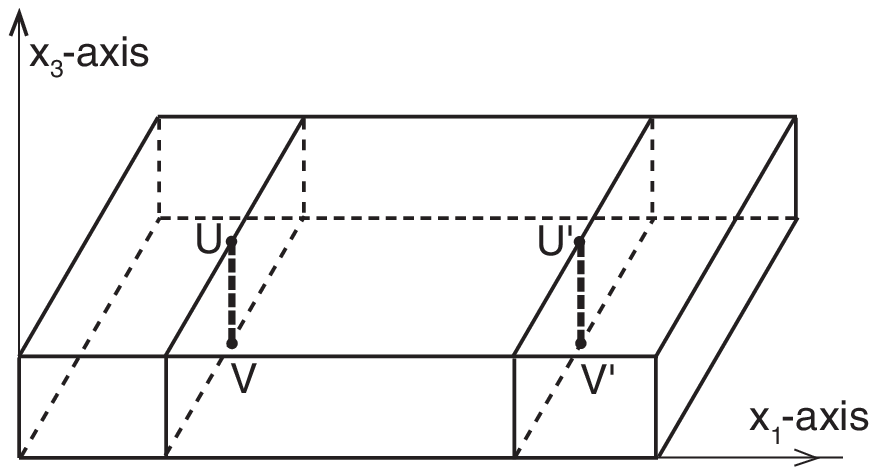}{4in} $$ \caption{The ball $B_1$}\label{X1} \end{figure}

We will present the ball $B_1$ as the rectangular parallelepiped,
see Figure \ref{X1}:

$$ B_1=\{ (x_1,x_2,x_3)\in \BR^3 \mid 0\le x_1 \le 6, 0\le x_2 \le
2, 0\le x_3\le 1\}.$$
\newcommand{\pr}{\mathrm{pr}}

 We suppose that  the knot intersects the interior of
$B_1$ in two vertical (i.e. parallel to the $x_3$-axis) straight
intervals $UV$ and $U'V'$, where $U=(1,1,1)$, $U'=(5,1,1)$, and
$V=(1,1,0), V'=(5,1,0)$. After an isotopy, we assume that the part
of $K$ outside the interior  of $B_1$ are 2 non-intersecting arcs
$\cc$ and $\cb$ on $S^2:=
\partial B_1$, where $\cc$ connects $U$ and $U'$, and $\cb$
connects $V$ and $V'$. Later we will describe explicitly the arc
$\cc$. If one cuts $S^2$ along the arc $\cc$, then one obtains a
disk, hence the other arc $\cb$ is uniquely determined by $\cc$,
up to isotopy.

\subsection{Skein module of complements of two-bridge knots}
Let $W$ be the top rectangle of $B_1$, $W =\{(x_1,x_2,x_3)\in X_1
\mid x_3=1\}$. Note that $X_1:= B_1\setminus (UV \cup U'V')$ is
the cylinder over a two-punctured disk $W\setminus \{U,U'\}$.
Hence $\cS(X_1)$ is isomorphic to the commutative algebra
$\R[x,x',y]$, as described in  subsection \ref{disk2}. Here $x$ is
a small loop circling $U$, $x'$ a small loop circling $U'$, and
$y=\partial W$. One of our main results is
\begin{theorem} The skein module $\cS(S^3\setminus \gb(p,m))$ is
free over $\cR$ with basis $\{x^ay^b, 0\le a, 0\le b\le
(p-1)/2\}$. \label{main2}
\end{theorem}

The remaining part of the section is devoted to a proof of this
theorem. Moreover, we will present  more explicit structures of
the skein module $\cS(S^3\setminus \gb(p,m))$.
\newcommand{\ve}{\varepsilon}

\subsection{Description of the two-bridge knot: the curve $\cc$} Let $\psi$ be the rotation by $180^o$ about the axis
$\{x_1=3, x_2=1\}$ (which is parallel to the $x_3$-axis and
passing through the center of the rectangle $W$). One has
$\psi(B_1)=B_1$.

For a set $Z\subset\BR^3$ let  $Z[a,b]$ be the part of $Z$ in the
strip $\{ a\le x_1 \le b\}$, i.e. $Z[a,b] := Z \cap
\{(x_1,x_2,x_3) \mid a\le x_1 \le b\}$. We will consider 3 pieces
$B_1[0,1], B_1[1,5]$, and $B_1[5,6]$ of $B_1$.

\begin{figure}[htpb] $$ \psdraw{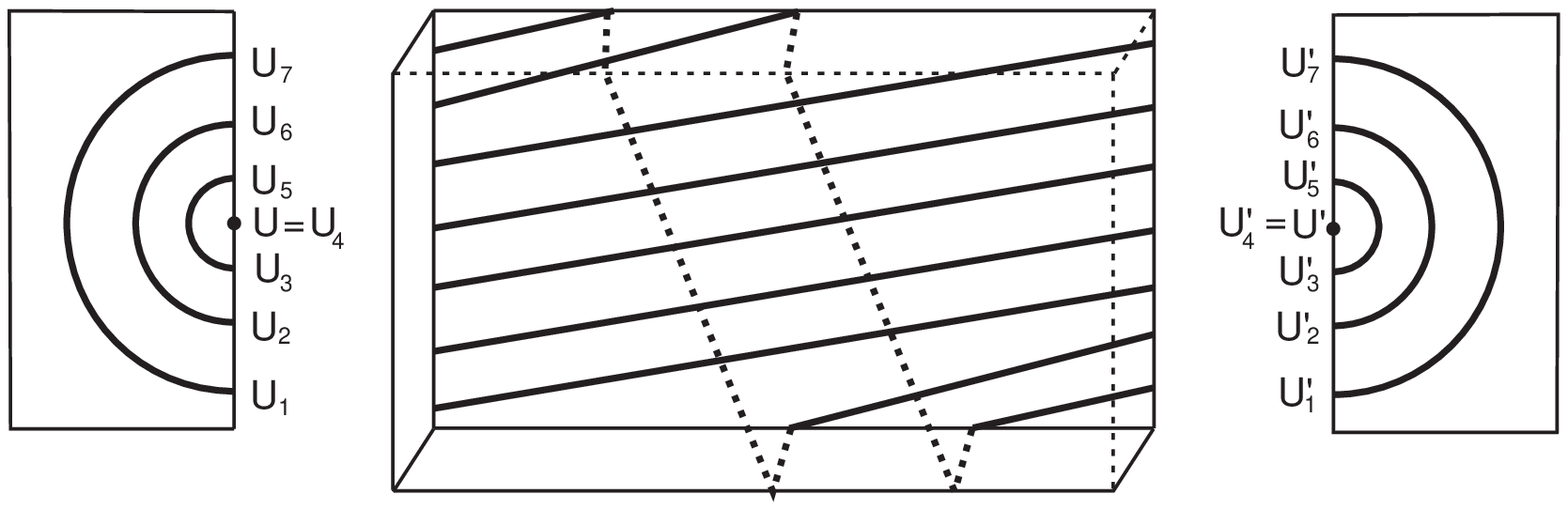}{5in} $$ \caption{$\cc_l, \cc_m$, and $\cc_r$ for $p=7,m=5$}\label{UL} \end{figure}

Recall that $W$ is the top rectangle of the rectangular
parallelepiped  $X_1$.
 On $W[0,1]$ let $\cc_l$ be the collection of
$(p-1)/2$ half-circles centered at $U$, with radii $\ve, 2\ve,
\dots (p-1)\ve/2$, where $\ve= 2/(p+1)$, see Figure \ref{UL}. The
end points of the half-circles, together with the center $U$, are
on a straight line; there are exactly $p$ of them, including $U$.
Mark them from bottom to top by $U_1,\dots, U_p$. (Thus
$U=U_{(p+1)/2}$.)

Similarly, on $W[5,6]$ let $\cc_r$ be the collection of $(p-1)/2$
half-circles centered at $U'$, with radii $\ve, 2\ve, \dots
(p-1)\ve/2$, where $\ve= 2/(p+1)$.  Mark the end-points of the
half-circles and $U'$ from bottom to top by $U'_1,\dots, U'_p$.
(Thus $U'=U'_{(p+1)/2}$.)

Recall that $S^2=\partial B_1$. On $S^2[1,5]$  we construct
$\cc_m$ as follows. Note that $S^2[1,5]$ is a cylinder over
$S^2[1,1]$: $S^2[1,5]= [1,5]\times S^2[1,1]$. Informally $\cc_m$
is the braid on the cylinder $S^2[1,5]= [1,5]\times S^2[1,1]$
representing the rotation by $(p-m)\pi /p$. Formally, we first
connect $U_1$ with $U'_{1+(p-m)/2}$ by a straight interval. Then
for all $i, 2\le i\le p$ connect $U_i$ with $U'_{i+(p-m)/2}$ by
non-intersecting arcs on $S^2[1,5]$; up to isotopy there is a
unique way to do so, see Figure ~\ref{UL}. There are in total $p$
arcs; denote them by $\cc_m$. We can assume that each arc in
$\cc_m$ always travel from left to right (no backwards traverse,
just like in the case of braids), and moreover, $\cc_m$ is
invariant under $\psi$.

\begin{figure}[htpb] $$ \psdraw{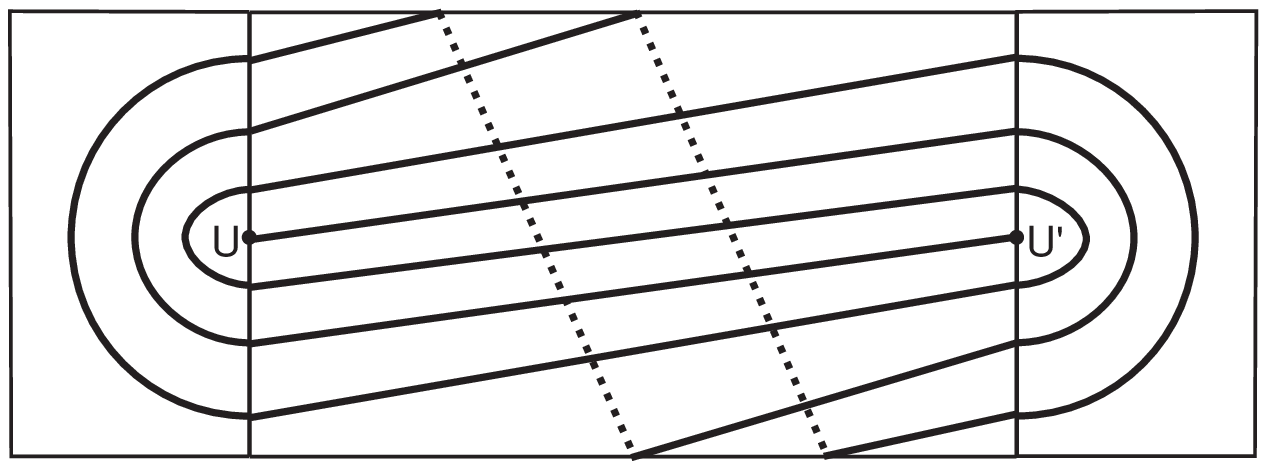}{4in} $$ \caption{The arc $\cc$}\label{ULL} \end{figure}

Let $\cc$ be the arc on $S^2$ obtained by combining $\cc_l, \cc_m$
and $\cc_r$; it  connects $U$ and $U'$ and is invariant under
$\psi$, see Figure \ref{ULL}. Up to isotopy there is a unique  arc
$\cb$ on $S^2$ connecting $V$ and $V'$.

\begin{proposition} The knot $K$ consisting of the arcs $\cc,\cb$, the
straight intervals $UV$ and $U'V'$ is the two-bridge knot
$\gb(p,m)$.
\end{proposition}
\newcommand{\mm}{\mathfrak m}
\newcommand{\lll}{\tilde {\mathfrak l}}
\begin{proof} It is clear that $K$ is a two-bridge knot. Let
$\tilde S^2$ be the 2-fold covering of $S^2=\partial B_1$ branched
along the 4 points $U,U',V,V'$. One can recognize the knot type by
looking at the homology class of the lift of $\cc$. Note that
$\tilde S^2$ is a torus, with the following preferred meridian and
longitude. The plane passing through $U,U',V,V'$ intersects
$S^2[0,1]$ in an arc $\mm$ that connect $U$ and $V$. The total
lift $\tilde \mm$ of $\mm$ is a closed curve on the the torus
$\tilde S^2$ which will serve as the meridian. The total lift
$\lll$ of the straight interval $UU'$ is another closed curve
serving as the longitude. It is is easy to see that $\tilde \mm$
and $\lll$ form a basis of $H_1(\tilde S^2,\BZ)$, and that the
total lift of the curve $\cc$, as a homology class, is equal to
$p\, \tilde\mm + (p-m)\lll$. According to \cite[Chapter
12]{Burde}, $K$ is a two-bridge knot of type $(p,m)$.
\end{proof}

\subsection{From $X_1$ to the knot complement $X$}
 Let $\ft$ be the boundary curve of a small normal
neighborhood of the arc $\cc$ in $S^2\equiv \partial B_1$. We can
assume that $\ft$ is invariant under $\psi$. Then $X = S^3
\setminus N(K)$ is obtained from $B_1$ by gluing a 2-handle to
along $\ft$.

\begin{figure}[htpb] $$ \psdraw{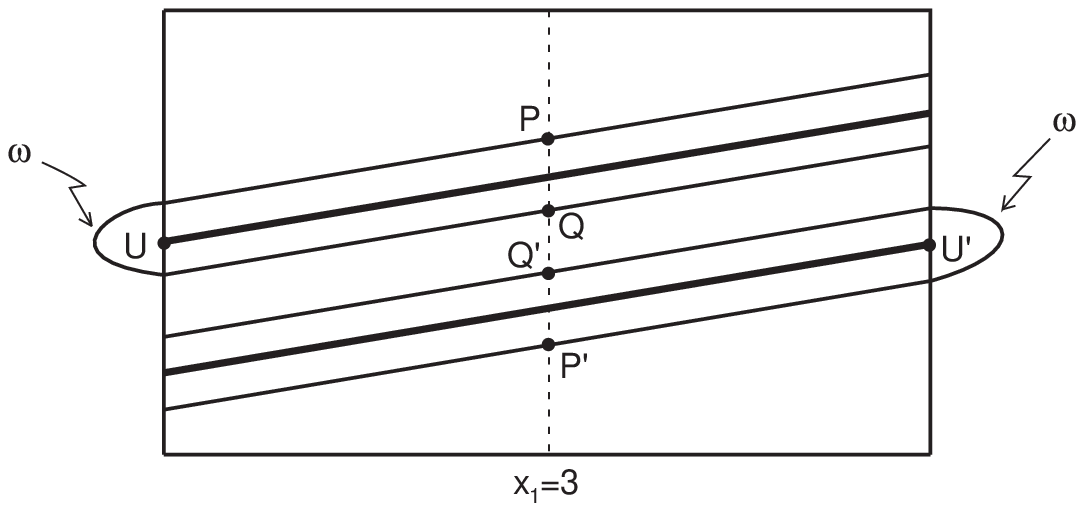}{4in} $$ \caption{The ball $X_1$}\label{ULLL} \end{figure}

Recall that on  $S^2[1,5]$ the curve $\cc$ consists of $p$ arcs,
and we assume that  the part of $\ft$ on $S^2[1,5]$ consists of
$2p$ arcs parallel and close to those of $\cc$. For example, the 2
arcs of $\cc$ on $S^2[1,5]$ containing $U$ and $U'$ are drawn in
Figure \ref{ULLL} by bold line,  and the part of $\ft$ near them
are drawn by lighter lines.  Let the plane $x_1=3$ intersect these
lighter lines (which are parts of $\ft$)  at the points
$P,Q,Q',P'$ (order from top to bottom), as in Figure \ref{ULLL}.
We have $\psi(P)=P', \psi(Q)=Q'$.

\subsection{Relative skein module} Let us recall the definition
of the relative skein module $\cS(X_1;P,Q')$ (see \cite{BullockL}).
A {\em type 1 tangle} is the disjoint union of a framed link and a
framed arc in $X_1$ such that  the parts of the arc near the two end
points are on the boundary $\partial X_1$, and the framing on these
parts are given by  vectors normal to $\partial X_1$. Type 1 tangles
are considered up to isotopy relative the endpoints. Then
$\cS(X_1;P,Q')$ is the $\R$-vector space generated by type 1 tangles
with endpoints at $P,Q'$ modulo the usual skein relations, like in
the definition of $\cS(X)$. One defines in a similar way the
relative Kauffman bracket skein module $ \cS(\Sigma;P,Q'):=
\cS(\Sigma\times [0,1];P,Q')$, where we identify $\Sigma\times
[0,1]$ with a collar of $\Sigma$ in $X_1$.

Since $\partial X_1=\Sigma$, there is a natural bilinear map
$\cS(\Sigma;P,Q') \otimes \cS(X_1) \to \cS(X_1;P,Q')$, where $\ell
\otimes \ell' \to \ell \star \ell'$, which  is the disjoint union
of $\ell$ and $\ell'$.

\newcommand{\ftt}{\mathfrak n}
\newcommand{\fc}{\mathfrak d}

The pair $P,Q$ divide $\ft$ into two arcs, the one that is fully
drawn in Figure \ref{ULLL} (and that goes around point point $U$)
is denoted by $\ftt_l$. Similarly, the pair $P',Q'$ divide $\ft$
into two arcs, the one that is fully drawn in Figure \ref{ULLL}
(and that goes around point point $U'$) is denoted by $\ftt_r$.
For $G,H$ among $P,Q,P',Q'$ let $\fc(GH)$ be the straight interval
connecting $G$ and $H$, whose interior is slightly pushed inside
the interior of $B_1$ (to avoid intersections with other arcs on
the boundary $\partial B_1$) and whose framing is given by vectors
normal to $\partial B_1$. Let $\fk_1$ be $\fc(P,Q')$; $\fk_2$ be
$\ftt_l$ followed by $\fc(QQ')$; $\fk_3$ be $\fc(PP')$ followed by
$\ftt_r$; and $\fk_4$ be $\ftt_l$ followed by $\fc(QP')$ then
followed by $\ftt_r$. In all cases the framing is given by vectors
normal to $\partial B_1$.

Using the skein relations one can simplify the arc part of
elements in $\cS(X_1;P,Q)$, showing that the arc part is one of
the four $\fk_i, i=1,2,3,4$. More precisely one has the following
lemma, which is a slightly different version of Lemma 3.1 in
\cite{BullockL}. The simple proof in \cite{BullockL} works for our
version as well.

\begin{lemma} The relative skein  $\cS(X_1;P,Q)$ is equal to the  union
$\displaystyle{\cup _{i=1}^4\left( \fk_i\star \cS(X_1)\right)}$.
\label{aa23}
\end{lemma}

\subsection{From $\cS(X_1)$ to $\cS(X)$ through sliding}

\subsubsection{} Let $\Sigma = S^2 \setminus \{U,U',V,V'\}$. Then
$\cS(\Sigma)$ is an algebra, and $\cS(X_1)=\cR[x,x',y]$ is a left
$\cS(\Sigma)$-module, since $\partial X_1= \Sigma$: For $\ell \in
\cS(\Sigma)$ and $\ell' \in \cS(X_1)$ the action $\ell\star \ell'
\in \cS(X_1)$ is the disjoint union of $\ell$ and $\ell'$. It is
easy to see that the action descends to an action of $\cS(\Sigma)$
on $\cR[x,y]\equiv \cR[x,x',y]/(x=x')$.
\subsubsection{Sliding}
Recall that  $X$ is obtained from $X_1$ by attaching a 2-handle
along the curve $\ft$. The embedding of $X_1$ into $X$ gives rise to
a linear map from $\cS(X_1)\equiv \cR[x,x',y]$ to $\cS(X)$. It is
known that the map is surjective, and its kernel $\cK$, see
\cite{Przytycki,BullockL}, can be described through slides as
follows.

Suppose $\fk$ is a type 1 tangle whose 2 endpoints are on $\ft$ such
that outside a small neighborhood of the 2 endpoints $\fk$ is in the
interior of $X_1$ and in a small neighborhood of the endpoints $\fk$
is on the boundary $\partial B_1=S^2$. The two end points of $\fk$
divide $\ft$ into 2 arcs $\ft_l$ and $\ft_r$. The loop $\ft$
partitions $S^2=\partial B_1$, which is a sphere, into 2 parts; the
one not containing $U,U'$ is called the {\em outside one}. Let us
isotope $\fk$ (relatively to the endpoints)  to $\fk'$ so that in a
small neighborhood of the endpoints  $\fk'$ is in the outside part
of $\ft$.

Let $sl(\fk)$ be $\fk'.\ft_l-\fk'.\ft_r$, considered as an element
of the skein module $\cS(X_1)$. Here $\fk'.\ft_l$ is the framed link
obtained by combining  $\fk'$ and $\ft_l$. Note that $sl(\fk)$ is
defined up to a factor $\pm t^{3n}, n\in \BZ$: The exchange $\ft_l
\leftrightarrow \ft_r$ changes the sign, and isotopies in
neighborhoods of the endpoints change the framing, which results in
a factor equal to a power of $(-t^3)$.

It is clear that as framed links in $X$, $\fk'.\ft_l$ is isotopic to
$\fk'.\ft_r$, since one is obtained from the other by sliding over
the 2-handle attached to the curve $\ft$. Hence we always have
$sl(\fk)\in\cK$. It was known that $\cK$ is spanned by all possible
$sl(\fk)$, where $\fk$ can be chosen among all type 1 tangles with
pre-given two endpoints on $\ft$.

\subsubsection{The kernel}

 Now we choose and fix the two end points of $\fk$: $P$ and $Q'$;
 i.e. $\fk\in \cS(X_1;P,Q')$. We will assume that $\ft_l$ contains $\ftt_l$. The kernel $\cK$ is spanned
 over $\cR$ by $sl(\fk), \fk\in \cS(X_1;P,Q')$. From the description of
$\cS(X_1;P,Q)$ in Lemma \ref{aa23}  we have

\begin{lemma} The kernel $\cK$ is equal to the $\cR$-span of
$\{sl(\fk_i)\star \cS(X_1), i=1,2,3,4 \}$.
\end{lemma}

We will call an element in $sl(\fk_i)\star \cS(X_1)$ a {\em
relation of type $i$}. Here $i=1,2,3,4$.

\begin{lemma}
\label{e3}
 One has
$$
\begin{aligned} sl(\fk_1) & = sl (\fc(PQ')) \\
sl(\fk_2) & = sl (\fc(QQ')) \\
sl(\fk_3) & = sl (\fc(PP')) \\
sl(\fk_4) & = sl (\fc(P'Q)) \\
\end{aligned}
$$
\end{lemma}
\begin{proof}
The first identity is a tautology. The last three follows from
trivially a simple isotopy of the links involved.
\end{proof}

\subsubsection{Simplifying the kernel} Let $\cK'$ be the $\cR$-span of
$sl(\fk_1)\star \cS(X_1) $ and $(x-x')\cS(X_1)$.

\begin{lemma}For every $\ell \in \cS(X_1)$ one has $\ell
-\psi(\ell) \in (x-x')\, \cS(X_1)$. \label{e4}
\end{lemma}
\begin{proof} Note that $\psi(x)=x',
\psi(x')=x$, and $\psi(y)=y$. Hence for any link $\ell$, the skein
$\psi(\ell)$, as an element in $\cS(X_1)=\R[x,x',y]$, is obtained
from the skein of $\ell$ by the involution $x\to x', x'\to x, y\to
y$. It follows that for any framed link $\ell$, and one has $\ell
-\psi \ell\in (x-x')\,\cS(X_1) \subset \cK'$.
\end{proof}

\newcommand{\dd}{\mathfrak l}
\begin{lemma} One has $\cK= \cK'$.
\end{lemma}
\begin{proof}
First we prove that $\cK'\subset \cK$. Since $sl(\fk_1)\star
\cS(X_1) $ is already in $\cK$, we need to show $(x-x')\,\cS(X_1)$
is in $\cK$. Note that $sl(\fc(PQ))$ is exactly $x-x'$. Moreover,
for any link $\ell$ in the interior of $X_1$ one has
$sl(\fc(PQ)*\ell)= (x-x')\ell$. Hence $(x-x')\, \ell\in \cK$, and
hence $\cK'\subset \cK$.

Now we prove $\cK\subset \cK'$.

By Lemma \ref{e3}  one has $sl(\fk_3)= sl(\fc(PP'))$. Since both
$\fc(PP')$ and $\ft$ is invariant under $\psi$, we have
$\psi(\fc(PP).\ft_l(PP')) = \fc(PP).\ft_r(PP')$, where
$\ft_l(PP')$ and $\ft_r(PP')$ are the two arcs of $\fc$ obtained
by dividing $\ft$ using the two points $P,P'$. Hence
$sl(\fc(PP'))\, \star \cS(X_1) = (\fc(PP).\ft_l(PP') -
\fc(PP).\ft_r(PP') )\, \star \cS(X_1)$ is in $(x-x')\, \cS(X_1)
\subset \cK'$ by Lemma \ref{e4}. The proof that all relations of
type 2 belongs to $\cK'$ is similar.

For type 4, by Lemma \ref{e3} we have $sl(\fk_4)=sl(\fc(P'Q))$.
Since $\psi(P)=P', \psi(Q)=Q'$, one has $ \psi(sl((\fc(P'Q))) = sl
(\fc(PQ'))=sl(\fk_1)$. Hence $[sl(\fk_4) -sl(\fk_1)]\star
\cS(X_1)$ belongs to $(x-x')\, \cS(X_1)$ by Lemma \ref{e4}. Thus
$sl(\fk_4) \star \cS(X_1)$ is in $\cK'$.
\end{proof}
 \subsection{Proof of theorem \ref{main2}} We have $\cS(X)=
 \cR[x,x',y]/\cK'$. Note that $\cR[x,x',y]/
 (x-x')\cS(X_1)=\cR[x,y]$.  Hence $\cS= \cR[x,y]/\cK''$, where $\cK''$ is the
 $\cR$-span of $sl(\fk_1)\star\cR[x,y]$. Note that there is a
 natural $\cR[x]$-module structure on $\cS(X)$: Here $x$ is a
 meridian, thus belongs to the boundary of $X$. Over $\cR[x]$,
 $\cR[x,y]$ is spanned by $1,y,y^2,\dots$. Hence $\cK''$, as a
 $\cR[x]$ module, is spanned by $ sl(\fk_1)\star y^k= (\fk_1.\ft_l-\fk_1.\ft_r)\star y^k ,
 k=0,1,2,\dots$.

Note that $\fk_1.\ft_r$ is the closure in the sense of Section
\ref{d3} of a braid on $(p+1)$ strands, while $\fk_1.\ft_l$ is the
closure  of a braid on $(p-1)$ strands. Moreover, $\fk_1.\ft_r
\star y^k$ is the closure of a braid on $(p+1) + 2k$ strands,
while $\fk_1.\ft_l \star y^k$ is the closure of of a braid on
$(p-1) + 2k$ strands. Lemma \ref{disk21} shows that
$(\fk_1.\ft_l-\fk_1.\ft_r)\star y^k$, as an element of $\cR[x,y]$,
has $y$-degree $(p+1)/2 + k$, with highest coefficient invertible
and of the form a power of $t$. Hence when we factor out
$\cR[x,y]$ by $\cK''$, we get a free $\cR[x]$-module with
representatives $y^l, l=0,1,2,\dots, (p-1)/2$ as a basis.
 This
completes the proof of Theorem \ref{main2}.

\begin{remark} The same proof shows that the theorem still holds true if we replace the
ground ring $\cR=\BC[t^{\pm1}]$ by $\cR_\BZ := \BZ[t^{\pm1}]$.
\end{remark}

\begin{corollary} For two-bridge knots one has $\epsilon(\cS(X))=
R(\cX(X))$, the ring of regular functions on the character
variety. \label{phu}
\end{corollary}

\begin{proof}By the result of \cite{Le1}, the ring $R(\cX(X))$ is
$\BC[\bar x, \bar y]/(\varphi(\bar x,\bar y))$, where
$\varphi(\bar x,\bar y)$ is a polynomial of $\bar y$-degree
$(p+1)/2$, with leading coefficient 1. Here $\bar x$, $\bar y$ are
respectively the traces of the loop $x,y$. The corollary follows
immediately.
\end{proof}

\section{Proof of Theorem \ref{main1}}
\label{proof1}

\subsection{The peripheral polynomial of two-bridge knots}Theorem
\ref{main2} about the structure of the skein module of complements
of two-bridge knots is used to prove the following.
\begin{proposition} For the two-bridge knot $K=\gb(p,m)$ the
peripheral polynomial $\beta_K$ is never 0 and has $L$-degree less
than or equal to $(p+1)/2$. \label{degree}
\end{proposition}

\begin{proof}By definition one has the following exact sequence of
$\cR[x]$-modules

\begin{equation} 0 \to \cP \hookrightarrow \cT^\sigma
\overset{\Theta}{\longrightarrow} \cS(X) \label{sao200}
\end{equation}

When $x=M+M^{-1}$, the field $\cR(M)$ of rational functions in $M$
 is a flat $\cR[x]$-module,  since
$\cR(M)$ contains the fractional field of $\cR[x]$ as a subfield.
Hence the following sequence, which is obtained from
(\ref{sao200}) by tensoring with $\cR(M)$, is exact

\begin{equation} 0 \to \cR(M)\otimes_{\cR[x]}\cP \, \hookrightarrow\,
\cR(M)\otimes_{\cR[x]}\cT^\sigma \,\overset{id \otimes
\Theta}{\longrightarrow}\cR(M)\otimes_{\cR[x]} \cS(X)
\label{sao201}
\end{equation}

Note that the first module $\cR(M)\otimes_{\cR[x]}\cP$ is exactly
$\tilde \cP$, the left-ideal extension of $\cP$ from $\cT^\sigma$
to $\tilde \cT$. It's easy to check that the second module
$\cR(M)\otimes_{\cR[x]}\cT^\sigma$ is $\tilde \cT$. One can now
rewrite (\ref{sao201}) as

\begin{equation} 0 \to  \tilde \cP \, \hookrightarrow\,
\tilde \cT \,\overset{id \otimes \Theta}{\longrightarrow}\,
\cR(M)\otimes_{\cR[x]} \cS(X) \label{sao202}
\end{equation}

The third module is a finite-dimensional $\cR(M)$-vector space; in
fact, its basis is $\{y^i, 0\le i\le \frac{p-1}{2}\}$, since
$\cS(X)$ is $\cR[x]$-free with the same basis, by Theorem
\ref{main2}. The middle module $\tilde \cT$ is an $\cR(M)$-vector
space of infinite dimension; in fact, its basis is $\{L^a, a\in
\BZ\}$. Thus the kernel $\tilde \cP$  is never 0, and hence its
generator $\b_K$ is not 0. Moreover the image of $(p+1)/2 +1$
elements $1, L, L^2, \dots, L^{(p+1)/2}$ are linearly dependent.
Hence there must be a non-trivial element in the kernel of
$L$-degree less than or equal to $(p+1)/2$.
\end{proof}

\begin{corollary} Every two-bridge knot $\bb (p,m)$  satisfies a
recurrence relation with $L$-degree less than or equal to
$(p+1)/2$.
\end{corollary}

This is because $\beta_K$ is divisible by the recurrence
polynomial $\alpha_K$. Note that the existence of recurrence
relations for {\em arbitrary} knots was established in \cite{GL}
by another method. But in \cite{GL} the $L$-degree is much larger.

\subsection{Conjecture 2 holds true for two-bridge knots}
\begin{proposition}\label{400}

Suppose for a knot $K$ the skein module $\cS(X)$ is free over
$\cR$ and $\epsilon(\cS(X)) = R(\cX(X))$. Then $\epsilon(\cP)=
\gp$, and hence $\epsilon(\beta_K)\overset{M}{=} B_K$.
\end{proposition}

\begin{proof} Consider again the exact sequence (\ref{sao200}),
but now as sequence of modules over $\cR=\BC[t^{\pm1}]$, a
principle ideal domain. By assumption, the last module $\cS(X)$ is
free over $\cR$. Hence when tensoring (\ref{sao200}) with any
$\cR$-module, one gets an exact sequence. In particular, tensoring
with $\BC$, considered as $\cR$-module by putting $t=-1$, one has
the exact sequence

$$ 0\to \epsilon(\cP) \hookrightarrow \epsilon(\cT^\sigma)
\overset{\epsilon(\Theta)}{\longrightarrow} \epsilon(\cS(X)).$$

Notice that $\epsilon(\cT^\sigma)= \gt^\sigma$, $\epsilon(\cS(X))=
R(\cX(X))$, and $\epsilon(\Theta) =\theta$. Thus $\gp$, being the
kernel of $\theta$, is equal to $\epsilon(\cP)$. The second
statement follows from Proposition \ref{t1}.
\end{proof}

From Theorem \ref{main2}, Corollary \ref{phu} and Proposition
\ref{400} we get

\begin{theorem} Conjecture \ref{weakconjecture} holds true for
two-bridge knots: $\epsilon(\cP)=\gp$ and $\epsilon(\beta_K)=
B_K$. \label{main3}
\end{theorem}
\subsection{Proof of Theorem \ref{main1}}

(a) is Proposition \ref{degree}.

(b)  One has $\epsilon(\beta_K)\overset{M}{=}B_K$ by Theorem
\ref{main3}. Thus the algebraic set
$\{\epsilon(\beta_K)=0\}\overset{M}{=} \{ B_K=0\}$ is
$M$-essentially equal to  $\{A'_K=0\}$, by Proposition \ref{500}.
Applying $\epsilon$ to (\ref{sao}) we get

\begin{equation} \epsilon(\beta_K) \overset {M}{=}\epsilon(\gamma) \,
\epsilon(\alpha_K). \label{300}
\end{equation}

which means  $\epsilon(\beta_K)$ is $M$-divisible by
$\epsilon(\alpha_K)$. Hence $\{\epsilon(\alpha_K)=0\}$ is
$M$-essentially an algebraic subset of
$\{\epsilon(\beta_K)=0\}\overset{M}{=} \{A'_K=0\}$.

(c) Suppose $A_K$ has $L$-degree $(p-1)/2$. Then $A'_K=(L-1)A_K$
 has $L$-degree $(p+1)/2$. By Proposition \ref{degree}, $\beta_K$ has
 $L$-degree less than or equal to $(p+1)/2$. But
 $\epsilon(\beta_K)$is $M$-essentially divisible by $A'_K$. Hence
 the $L$-degree of $\beta_K$ must be exactly $(p+1)/2$, and also

$$ \epsilon(\beta_K) \overset{M}{=} (L-1) A_K.$$

Combining with (\ref{300}) we have

\begin{equation} A_K \overset {M}{=}\epsilon(\gamma) \,
\frac{\epsilon(\alpha_K)}{L-1}. \label{600}
\end{equation}

Recall that $\frac{\epsilon(\alpha_K)}{L-1}$ is a polynomial by
Proposition \ref{divi}.  From (\ref{sao}) and the fact that the
$L$-degree of $\alpha_K$ is bigger than 1 (Proposition \ref{deg2})
it follows that the $L$-degree of $\gamma$ is less than $(p-1)/2$,
which is the $L$-degree of $A_K$. Hence if $A_K$ is
$\BZ$-irreducible, from (\ref{600}) one must have $\epsilon
(\gamma)\overset{M}{=}1$ and
$\frac{\epsilon(\alpha_K)}{L-1}\overset{M}{=} A_K$. This completes
the proof of Theorem \ref{main1}.


\begin{thebibliography}{10000}

\bibitem[BG]{BG}D. Bar-Natan and S. Gaorufalidis,
{\em On the Melvin-Morton-Rozansky conjecture}, Invent. Math. {\bf
125} (1996), no. 1, 103--133.

\bibitem[BZ]{Burde}
G. Burde and H. Zieschang, {\em Knots}, de Gruyter Studies in
Mathematics, 5. Walter de Gruyter \& Co., Berlin, 2003.


\bibitem[Bu1] {Bullock} D. Bullock,
{\em The $(2,\infty)$-skein module of the complement of a $(2,2p+1)$
torus knot},  J. Knot Theory Ramifications {\bf 4} (1995), no. 4,
619--632.

 \bibitem[Bu2] {Bullock4} D. Bullock, {\em Rings of SL2C -characters and the Kauffman bracket skein module},
 Comment. Math. Helv. {\bf 72} (1997), no. 4, 521--542.

\bibitem[BL] {BullockL} D. Bullock and W. Lo Faro,
{\em The Kauffman bracket skein module of a twist knot exterior},
 Algebr. Geom. Topol.  {\bf 5} (2005), 107--118 (electronic).


\bibitem[CCGLS] {CCGLS} D. Cooper, M. Culler, H. Gillett, D.D. Long and P.B. Shalen,
{\em Plane curves associated to character varieties of 3-manifolds},
Inventiones Math. {\bf 118} (1994), pp. 47--84.

\bibitem[CL] {CL}
D. Cooper and D. Long, Representation theory and the $A$- polynomial
of a knot, Chaos, Solitons and Fractals, {\bf 9} (1998) no 4/5,
749--763.

\bibitem[CS] {Culler}  M. Culler and P.B. Shalen, {\em Varieties of group representations and splittings of $3$-manifolds},
 Ann. of Math. (2) {\bf 117} (1983), no. 1, 109--146.


\bibitem[FG]{FG}C. Frohman and R. Gelca, {\em
Skein modules and the noncommutative torus}, Trans. Amer. Math. Soc.
{\bf 352} (2000), no. 10, 4877--4888.

\bibitem[FGL]{FGL}C. Frohman, R. Gelca, and W. Lofaro,
{ \em The A-polynomial from the noncommutative viewpoint},  Trans.
Amer. Math. Soc. {\bf 354} (2002), no. 2, 735--747.

\bibitem[Ga]{Ga2} S. Garoufalidis, {\em
On the characteristic and deformation varieties of a knot},
 Proceedings of the Casson Fest,  291--309 (electronic), Geom. Topol. Monogr., 7, Geom. Topol. Publ., Coventry, 2004.

\bibitem[GL]{GL} S. Garoufalidis and T. T. Q. L\^e,
{\em The colored Jones function is q-holonomic},  Geom. Topol. {\bf
9} (2005), 1253--1293 (electronic).

\bibitem[Ge]{Gelca} R. Gelca,
{\em On the relation between the $A$-polynomial and the Jones
polynomial}, Proc. Amer. Math. Soc. {\bf 130} (2002), no. 4,
1235--1241 (electronic).

\bibitem[HS]{Hoste} J.
Hoste and P. Shanahan, {\em
 Commensurability classes of twist knots},  J. Knot Theory Ramifications  {\bf 14}  (2005),  no. 1, 91--100.

\bibitem[Hi]{Hikami}K. Hikami,
{\em Difference equation of the colored Jones polynomial for torus
knot},  Internat. J. Math.  {\bf 15}  (2004),  no. 9, 959--965.

\bibitem[Jo]{Jones}
 V.F.R. Jones, {\em Polynomial invariants of knots via von Neumann
 algebras}, Bull. Amer. Math. Soc., {\bf 12} (1985), 103--111.

\bibitem[Ka]{Kauffman}
L. Kauffman, {\em States models and the Jones polynomial}, Topology,
{\bf 26} (1987), 395--407.

\bibitem[Le1]{Le1} T. T. Q. L\^e,  {\em Varieties of representations and their
subvarieties of cohomology jumps for knot groups}, (Russian)  Mat.
Sb.  {\bf 184}  (1993), no. 2, 57--82;  translation in  Russian
Acad. Sci. Sb. Math.  {\bf 78}  (1994),  no. 1, 187--209.
\bibitem[Le2]{Le10}T. T. Q. L\^e, {\em On the skein modules and the A-polynomial of two-bridge
knots}, to appear.

  \bibitem[LM]{Lubotzky}
A. Lubotzky, A. Magid, {\em Varieties of representations of finitely
generated groups}, Memoirs of the AMS {\bf 336} (1985).

  \bibitem[Li]{Lickorish}W. B. R. Lickorish, {\em An Introduction to Knot Theory}, Springer, GTM {\bf 175}, 1997.


\bibitem[MM]{MM}
P. M. Melvin and H. R. Morton,  {\em The coloured Jones function},
Comm. Math. Phys. {\bf 169} (1995), no. 3, 501--520.

\bibitem[Pr]{Przytycki}
J. Przytycki, {\em Fundamentals of Kauffman bracket skein modules},
Kobe J. Math. {\bf 16} (1999), 45--66.

\bibitem[PS]{PS}
J. Przytycki, A. S. Sikora,  {\em On skein algebras and ${\rm Sl}\sb
2( C)$-character varieties}, Topology {\bf 39} (2000), 115--148.

\bibitem[Ta]{Takata} T. Takata, {\em
The colored Jones polynomial and the A-polynomial for twist
knots}, preprint math.GT/0401068.

\bibitem[Tu]{Tu}  V.~G. Turaev, {\it Quantum invariants of
knots and 3-manifolds}, de Gruyter Studies in Mathematics 18, Walter
de Gruyter, Berlin New York 1994.

\end{thebibliography}
\end{document}